\documentclass[12pt]{article}

\title{Some Remarks on Some Strongly Coupled Reaction-Diffusion
Equations}
\author{Toka Diagana}
\date{May 8, 2003}
\begin{document}
\maketitle
%\tableofcontents

\noindent Abstract:{\small The primary goal of this paper is
to characterize solutions to coupled reaction-diffusion
systems. Indeed, we use operators theory to show that under suitable assumptions, then the
system given by $$ \; u_{t} = M \; \Delta \; u + F(u)$$
have solutions. As applications, we consider a mathematical
model arising in Biology and in Chemistry.}

\bigskip


\vspace{0.2cm}


\section{Introduction}
The primary goal of this paper is to make an investigation on a
particular type of partial differential systems, that is, the
strongly coupled reaction-diffusion system. Recall that early
investigations on this problem are due to Fichera, see, e.g., [8],
and recently by Amann, see, e.g., [1, 2, 3]. Also, recall that in
most of publications on the proposed problem, the diffusion matrix
$M$ is supposed to be diagonal and sometimes with positive
entries. Obviously such an assumption cannot be applied to some
interesting cases arising in several fields such as in
Biology, Chemistry, and Ecology. In this paper, we consider the general
case, that is, we assume that $M$ is any matrix without any
restrictions on its entries. However it will be shown that if all
eigenvalues of $M$ belong to $S = \{ z \in {\mathbf{C}}: \; \Re e
z \geq 0 \}$, then the reaction-diffusion systems have a unique solution. 
The main idea of our investigation is based on 
operators theory, especially unbounded normal operators and related
semi-groups of contraction. The strongly coupled reaction-diffusion system
is defined as
\begin{equation}
u_{t} = M \; \Delta \; u + F(u),
\end{equation}
where $M$ is a $d \times d$ real matrix and $F : {\mathbf{R}}^d
\mapsto {\mathbf{R}}^d$ is of class $C^2$. Throughout this paper we will
assume that $d$ is even; the case where $d$ is odd will be
investigated elsewhere. However, the author expects to use
a similar method as in this paper.  

As stated above, we propose to solve
Eq.(1) using unbounded normal operator method. Under appropriate hypotheses, we prove that Eq.(1) has a unique
solution. 

Let us consider the needed background for it. Let
$A$(not necessarily bounded) be a normal operator in the (complex)
space Hilbert ${\mathcal{H}}$ ( $A$ is densely defined, and closed such
that $A A^* = A^* A$). Using the spectral theorem for unbounded
normal operators, see, e.g., [4, or 11, pp. 348-355], it is
well-known that $A = A_{1} - i A_{2}$, where $A_{1}$ and $A_{2}$,
are respectively the real part and minus the imaginary part of $A$.  Notice
that $A_{k}$'s are respectively self-adjoint
operators. It is also well-known that if we assume
that $A_{k}$'s are nonnegative self-adjoint operators, then $i A$
is an m-accretive operator. Thus $(i A)$ is the infinitesimal
generator of a contraction semi-group, see, e.g., [4, or 10
Corollary 4.4, p. 15]. To examine Eq.(1), we study the
linear part of it, that it, the diffusion operator $L = M \Delta$
in the Hilbert space $H = [ L^2 ( \Omega) ]^d$ where $\Omega$ is
a bounded subset of ${\mathbf{R}}^p$ with a smooth boundary. In
fact, we need to study the following problem

\[ (S_{0}) \; \; \; \left \{ \begin{array}{ll}
u_{t} = - M \Delta u \\ 
u_{|_{\partial{\Omega}}} = 0
\end{array}
\right. \]

More generally, let
$A$, $B$, $C$, and $E$ be unbounded normal operators in $H$, and
consider the following system

\[ (S_{1}) \; \; \; \left
\{ \begin{array}{llll}
u_{t} = A u + B v \\ 
v_{t} = C u + E v \\ 
u(0,x) = u_{0}(x) \\
v(0,x) = v_{0}(x)
\end{array} 
\right. \]

Such a system is equivalent to the following:

\[ (S_{2}) \; \; \; \left \{
\begin{array}{ll}
X_{t} = T \; X \; \; \mbox{with} \; \; X = (u,v) \\ 
X(0,x) = X_{0}(x) = (u_{0}(x) , v_{0}(x))
\end{array}
\right. \]

where the unbounded
matrix operator $T$
is defined by
$$ T \; : [D(A) \cap D(C)] \oplus [D(B) \cap D(E)]
\rightarrow H_{0} \oplus H_{0}$$ and

$$
T =  \; \left(
\begin{array}{cc}
A  &  B \\
C &  E
\end{array}
\right)
$$

Set $ T \; = \; - i \; \tilde{T} $. Thus, the unbounded
operator $\tilde{T}$ is defined by
$$
\tilde{T} =  \; \left(
\begin{array}{cc}
i \; A  &  i \; B \\
i \; C &  i \; E
\end{array}
\right)
$$
It follows that $(S_{1}) \Leftrightarrow (S_{2}) \Leftrightarrow
(S_{3})$, where $(S_{3})$ is given by
\[ (S_{3}) \; \; \; \left \{
\begin{array}{ll}
X_{t} = - i \tilde{T} \; X \; \; \mbox{with} \; \; X = (u,v) \\ 
X(0,x) = X_{0}(x) = (u_{0}(x) , v_{0}(x))
\end{array}
\right. \]

The next step is to show that $\tilde{T}$ is a normal operator, which implies that $- T$ is an m-accretive operator.

\section{Diffusion Equation}
In this section, we show that the general problem $(S_{1})$
admits a unique solution under appropriate hypotheses on $A$, $B$,
$C$, and $E$. As particular case, the equation $(S_{0})$ will be
considered.

\bigskip

{\textbf{Theorem 2.1}}.- {\it Let $A$, $B$, $C$, and $E$ be
unbounded normal operator in the Hilbert space $H$. Assume the
following conditions hold true
\begin{enumerate}
\item $N(A) = \{0\}$, and $N(E) = \{0\}$
\item $D(A) \cap D(C)$ and $D(B) \cap D(E)$ are dense in $H$
\item $A^{-1} B$ and $C^{-1} E$ are closed operators in $H$
\item $A$, $B$, $C$, and $E$ commute each other
\end{enumerate}
where $N(A)$ and $N(E)$ are respectively the Kernels of the
operators $A$ and $E$.
Then the matrix operator $\tilde{T}$ is normal in $H \oplus H$.}

\bigskip

Proof. Recall that the matrix operator $\tilde{T}$ is defined in
$H \oplus H$ by, $D( \tilde{T}) = [D(A) \cap D(B) ] \oplus [ D(C)
\cap D(E)], \; \; \tilde{T} (u,v) = < i(Au + Bv) , i(C u + E v)
> \; \; \forall (u,v) \in D(\tilde{T})$. Therefore, it is a densely
defined operator in $H \oplus H$, according to 
assumption (2). Let $(u_{n} , v_{n})$ be a sequence in $D(
\tilde{T})$ such that $(u_{n} , v_{n})$ converges to $(u,v)$ and
$\tilde{T} (u_{n} , v_{n})$ converges to $(i \xi , i \eta)$ in $H
\oplus H$. In other words, $A u_{n} + Bv_{n}$ and $ C u_{n} + E
v_{n}$ converge to $ \xi$ and $\eta$ respectively. Since the
kernel $N(A) = \{0\}$ (according to (1) ), then $A^{-1} B v_{n}
\longrightarrow A^{-1} \xi - u$. Now, since $A^{-1} B$ is closed
then $v
\in D(A^{-1} B) = D(B)$ and $ A^{-1} B v = A^{-1} \xi - u$. In
addition $u \in D(A)$, since $u = A^{-1} \xi - A^{-1} Bv \in
D(A)$. Thus, it easily follows that $(u,v) \in D(A) \oplus D(B)$ and $A u
+ Bv = \xi$. Using a similar argument yields $(u,v) \in D(C)
\oplus D(E)$ and $Cu + E v = \eta$. Therefore, $\tilde{T}$ is a
closed operator. Now, we have $$ \tilde{T} {\tilde{T}}^* =  \;
\left(
\begin{array}{cc}
A^* A + C^* B  &  B^* A + E^* B \\
A^* C + C^* E & B^* C + E^* E
\end{array}
\right) $$ and,  $$ {\tilde{T}}^* \tilde{T} =  \; \left(
\begin{array}{cc}
AA^* + C B^*  &  B A^* + E B^* \\
A C^* + C E^* & B C^* + E E^*
\end{array}
\right) $$ Clearly $A A^* = A^* A$, and $E E^* = E^* E$( $A$ and $E$
are normal operators). Let show that $C^* B =
B^* C$. A similar argument can be used to show that $B A^* = A
B^*$, $E^* B = E B^*$, and $C^* E = C E^*$. Let us write $C =
C_{1} - i C_{2}$ and $B = B_{1} - B_{2}$ as stated in the
introduction of this paper. Thus $C^* = C_{1} + i C_{2}$ and
$B^* = B_{1} + i B_{2}$. Therefore $C^* B = (C_{1} B_{1} + C_{2}
B_{2}) + i ( C_{2} B_{1} - C_{1} B_{2})$ ; in the same way $B^* C
= ( B_{1} C_{1} + B_{2} C_{2}) + i ( B_{2} C_{1} - B_{1} C_{2})$.
Now, since $B C = C B$ (according to (4)), we have $C_{p} B_{q} =
B_{q} C_{p}$ for $p,q = 1, 2$. It follows that $C^* B = B^* C$.
In summary, we have $ {\tilde{T}}^* \tilde{T} = \tilde{T}
{\tilde{T}}^*$. The proof is complete.

\bigskip

{\textbf{Corollary 2.2}} Under previous assumptions. The operator
$\tilde{T}$ can be decomposed as $ \tilde{T} = T_{1} - i T_{2}$
where $T_{1}$ and $T_{2}$ are respectively self-adjoint operators. Assume that
both $T_{1}$, $T_{2}$ are nonnegative operators. Then the operator $- T
= i \tilde{T}$ is m-accretive. In addition $S_{1}$ admits a unique
solution.

\bigskip

Prof. Since $\tilde{T}$ is a normal operator and that both its
real and minus imaginary parts are nonnegative self-adjoint
operators, then $i \tilde{T} = - T$ is m-accretive, see, eg., [10,
Corollary 4.4, p. 15]. It follows that the system $(S_{1})$
admits a unique solution. Equivalently both $(S_{2})$ and $(S_{2})$ admit unique
solutions.

We apply previous results to the problem $(S_{0})$.
Assume that $d = 2n$ and set
\begin{equation}
A = - M_{1} \Delta, \; \; B = - M_{2}
\Delta
\end{equation}
$$
C = - M_{3} \Delta, \; \; \mbox{and} \; \; E = -
M_{4} \Delta
$$
where $M_{k}$ ($k = 1, 2, 3, 4$) is a $n \times
n$-matrix. Consider the following problem
\[(S_{4}) \; \;  \left \{
\begin{array}{ll}
X_{t} = - M \Delta \; X \; \; \mbox{with} \; \; X = (u,v) \\ 
X(0,x) = X_{0}(x) = (u_{0}(x) , v_{0}(x))
\end{array}
\right.
\]
where
$ M \Delta \; : [H^2(\Omega) \cap H_{0}^{1}(\Omega)]^{2n}
\rightarrow [ L^2(\Omega)]^{2n}$ is defined as

$$
M \Delta = \; \left(
\begin{array}{cc}
M_{1} \Delta  &  M_{2} \Delta \\
M_{3} \Delta &  M_{4} \Delta
\end{array}
\right)
$$
Set $T = - M \Delta$ and $- T = i \tilde{T}$. Therefore
$\tilde{T}$ is defined as
$$
\tilde{T} = \; \left(
\begin{array}{cc}
-i M_{1} \Delta  &  -i M_{2} \Delta \\
-i M_{3} \Delta &  -i M_{4} \Delta
\end{array}
\right)
$$
According to theorem 2.1, $\tilde{T}$ is a normal
operator. We have the following result.

\bigskip

{\textbf{Theorem 2.3}} {\it Under previous assumptions. Assume
that $M$ is not the zero matrix. In addition if all eigenvalues of
$M$ belong to the set $S = \{ z \in {\mathbf{C}} : \; \Re e z \geq
0 \}$. Then the problem $(S_{4})$ admits a unique solution.}

\bigskip

Proof. The operators $A$, $B$, $C$, and $E$ given in Eq.(2)
satisfy assumptions (1)--(2)--(3)--(4) of the theorem 2.2.
It turns out
that $\tilde{T}$ is a normal operator. Now let us show that if all
eigenvalues of $M$ belong to $S = \{ z \in {\mathbf{C}}: \; \Re e
z \geq 0 \}$, then $i \tilde{T} = M \Delta$ is m-accretive. Let
$\lambda_{1}, \; \lambda_{1}, ..., \lambda_{r}$ ( $k_{1} + k_{2} +
...+ k_{r} = 2n$) be eigenvalues of $M$. Following the Jordan
decomposition method for $M$, it is well-known that $M \Delta $
can be decomposed as
$$
M \Delta = \; \Pi \; \left(
\begin{array}{cccc}
J_{k_{1}}(\lambda_{1}) \Delta  & 0 & ... & 0\\
0 & J_{k_{2}}(\lambda_{2}) \Delta & ... & 0 \\
... & ... & ... & ... \\
0 & 0 & ... &  J_{k_{r}}(\lambda_{r}) \Delta
\end{array}
\right) \Pi^{-1}
$$
where $\Pi$ is a nonsingular matrix, and

$$
J_{k}(\lambda) \Delta \; = \; \left(
\begin{array}{cccc}
\lambda \Delta  & \Delta & ... & 0\\
0 & \lambda \Delta & \Delta & .. \\
... & ... & ... & ... \\
0 & 0 & ... &  \lambda \Delta
\end{array}
\right), \; k \geq 2; \; J_{1}(\lambda) = [\lambda \Delta], \; k =
1
$$
Thus, all eigenvalues of $M$ belong to $S = \{ z \in \mathbf{C} :
\; \Re e z \geq 0\}$ if and only if $i \tilde{T} = M \Delta$ is
m-accretive. Therefore, if all eigenvalues of $M$ belong to $S$
then $ i \tilde{T} = M \Delta$ generates a contraction semi-group.
In such a case $(S_{4})$ admits a unique solution. Equivalently
$(S_{0})$ admits a unique solution, since $(S_{0})$ is a
particular case of $(S_{4})$.

\section{Coupling Problem}
Assume $d = 2n$ and consider the coupling of the diffusion 
with the reaction term, that is, the system given by Eq.(1). We
define the following operators

\[ \left \{
\begin{array}{ll}
D(T) =  [H_{0}^{1}(\Omega) \cap H^2(\Omega)]^d \} \\ 
T u =  M \Delta u \; \; \; \forall u \in D(T)
\end{array}
\right. \]

where $M$ is the $d \times d$-matrix given in the Eq.(1).

\[ \left \{
\begin{array}{ll}
D(R) = \{ u \in [L^2(\Omega) ]^d : \; \; F(u) \in [L^2(\Omega) ]^d \} \\ Ru(t,x) =  F(u(t,x)) \; \; \; \mbox{a.e} \; \; u \in D(R)
\end{array}
\right. \]

We will make the following hypotheses

\bigskip

\noindent $(H_{0}) \; \; \; \mbox{all eigenvalues of the matrix
M belong to} \; \; S = \{ z \in {\mathbf{C}} : \; \Re e z \geq 0
\}$

\bigskip

\noindent $(H_{1}) \; \; \; \mbox{Assume the operator R is
m-accretive, and that} \; \; 0 \in R(0)$

We have the following.

\bigskip

{\bf Theorem 3.1}. {\it Under assumptions $(H_{0})$ and $(H_{0})$.
Then the Eq.(1) admits a unique solution.}

\bigskip

Proof. The main idea is to show that the nonlinear operator given
by $T + R$ is m-accretive in $H = [L^2(\Omega) ]^d$. Consider
Yosida's approximation for $R$. It is defined as
\begin{equation}
R_{\lambda} = \displaystyle{\frac{1}{\lambda}} [ I - ( I + \lambda
R)^{-1}], \; \; \; \lambda \; > \; 0
\end{equation}
It is well-known that $R_{\lambda}$ is m-accretive and that
$\displaystyle{\frac{1}{\lambda}}$-Lipschitz in $H$. Now consider
the following equation
\begin{equation}
\varepsilon u_{\lambda} + T u_{\lambda} + R_{\lambda} u_{\lambda}
= v, \; \; \; \mbox{and} \; \; \; u_{\lambda} = 0, \; \; \mbox{on}
\; \; \partial \Omega
\end{equation}
Since $T + R_{\lambda}$ is m-accretive, see, e.g., [5], then
Eq.(4) admits a unique solution $u_{\lambda} \in D(T)$ for any $v
\in H$, and $\lambda > 0$. We also know the family
$(u_{\lambda})_{\lambda
> 0}$ is bounded by $\displaystyle{\frac{1}{\varepsilon}} \; \| v
\|_{H}$. Using the fact $R_{\lambda}$ is m-accretive, $R_{\lambda}
0 = 0$, and by integration by parts it easily follows that
\begin{equation}
\int_{{\Omega}^d} Tu R_{\lambda} u dx \geq 0, \; \; \; \forall u
\in [H_{0}^{1}(\Omega) \cap H^2(\Omega)]^d
\end{equation}
Thus, multiplying Eq.(4) by $T u_{\lambda}$, and from Eq.(5), it
turns out that $(T u_{\lambda})$, and $(R_{\lambda} u_{\lambda})$
are bounded. From the compactness embedding, $[ H^2(\Omega)]^d
\hookrightarrow [L^2(\Omega)]^d$, and the fact $(u_{\lambda})$,
$(T u_{\lambda})$, and $(R_{\lambda} u_{\lambda})$ are bounded, it
turns out that: $u_{\lambda}$ strongly converges to $u$, $(T
u_{\lambda})$ weakly converges to $\xi$, and $(R_{\lambda}
u_{\lambda})$ weakly converges to $\eta$, as $\lambda$ approaches
to $0$ in $H$. Since $T$ is closed, then $Tu = \xi$. Since $R$ is
m-accretive, then $R u = \eta$, see, e.g., [5]. In summary $T+R$
is m-accretive under assumptions $(H_{0})$ and $(H_{0})$.
Therefore the algebraic sum(see, e.g., [6]) $(T+R)$ generates a nonlinear
contraction semi-group, that is the Eq.(1) admits a unique
solution.

\section{Applications}
In this section, we consider a model considered in [9]. The
problem we will study represents a mathematical model
describing various chemical and biological phenomena. In [9],
Lyapunov functionals have used to prove a global existence of
unique solutions. Here, we use the method described above to prove
that the given problem admits a unique solution, under suitable
assumptions.

Our model is described as
\[ (M) \; \;
\left \{
\begin{array}{llll}
\noindent u_{t}  - \alpha \Delta u - \beta \Delta v = - \sigma
f(u,v)
\; \; \mbox{in} \; \; (0, \infty) \times \Omega \\
v_{t}  - \gamma \Delta u - \alpha \Delta v = \rho f(u,v) \; \; \mbox{in} \; \; (0, \infty) \times \Omega \\
u_{\nu} = v_{\nu} = 0 \; \; \mbox{on} \; \; (0, \infty)
\times \partial \Omega \\
u(0,x) = u_{0}(x), \; \; \; \; v(0,x) = v_{0}(x) \; \; \mbox{in}
\; \; \Omega
\end{array}
\right.
\]
where $u(t)$ and $v(t)$ represent either chemical concentrations
or biological population densities, $\Omega$ is a bounded open
subset of class $C^1$ in ${\mathbf{R}}^n$, $u_{\nu}$ (respectively
$v_{\nu}$) denotes the outward normal derivative on $\partial
\Omega$, and $\alpha, \; \beta, \gamma, \; \rho$, and $\sigma$ are
positive constants . In [9], the following hypothesis is made

\begin{equation}
2 \alpha > ( \beta + \gamma)
\end{equation}
The (M) can be expressed as
\[ (N) \; \;
\left \{
\begin{array}{ll}
(u_{t} , v_{t}) = M \Delta (u,v) + F (u,v) \\
(u_{\nu} , v_{\nu}) = (0,0) \\
(u (0,x), v(0,x)) = (u_{0}(x) , v_{0}(x))
\end{array}
\right.
\]
where $F(u,v) = ( - \sigma f(u,v)\; , \; \rho f(u,v))$, and
$$
M =  \; \left[
\begin{array}{cc}
\alpha  &  \beta\\
\gamma &  \alpha
\end{array}
\right]
$$
It is obvious to see that all eigenvalues of of the diffusion
matrix $M$ are given as, $EV(M) = \{ \alpha + \sqrt{\beta \gamma}
\; , \; \alpha - \sqrt{\beta \gamma} \}$. Since all eigenvalues of
$M$ are nonnegative, then
\begin{equation}
\alpha \; > \sqrt { \beta \; \gamma}
\end{equation}

Clearly Eq.(6) implies Eq.(7). Indeed,
$\displaystyle{\frac{1}{2}} ( \beta + \gamma) \geq \sqrt{\beta \;
\gamma}$. Therefore, instead of considering Eq.(7), we will only
assume that Eq.(6) holds.

Consider the Hilbert space $H = [L^2(\Omega)]^2$ and set

$$D(T) = [ H_{0}^{1}(\Omega) \cap H^{2}(\Omega)
]^2$$ and $$T (u,v) = M (\Delta u , \Delta v)$$In the same way,
define
$$D(R) = \{ (u,v) \in [L^2(\Omega) ]^2 : \; \;
(-\sigma f(u,v) , \rho f(u,v)) \in [L^2(\Omega) ]^2 \}$$ and
$$
R(u, v) = ( - \sigma f(u(t,x), v(t,x)) \; ,\; \rho f( u(t,x) ,
v(t,x)) \; \; \; \mbox{a.e} \; \; u, v \in D(R)$$

We will make the following assumption
\begin{equation}
f(0,0) = 0
\end{equation}
For instance from the fact that $R$ is accretive, the following holds
\begin{equation}
- \sigma u f(u,v) + \rho v f(u,v) \geq 0, \; \; \forall (u,v) \in
D(R)
\end{equation}
More generally, assume that $f$ is given such that $R$ is a
nonlinear m-accretive operator in $[L^2(\Omega)]^2$. Thus, we have
the following.

\bigskip

{\textbf{Proposition 4.1}}. Under Eq.(6), and Eq.(8), then the
problem described in (M) admits a unique solution.

\bigskip

Proof. Obvious as consequences of Theorem 2.3 and Theorem 3.1.

\vspace{0.8cm}

\begin{center}
{\Large References}
\end{center}

\vspace{0.5cm}

\begin{enumerate}

\item H. Amann, {\it Global Existence for Semilinear Parabolic
Systems}, J. Reine Angew. Math {\bf 360}, 47-83 (1985).

\item H. Amann, {\it Dynamic Thoery of Quasilinear Oarabolic
Systems. III Global Existence}. Math. Z. {\bf 202}, 2, 219-250
(1989).

\item H. Amann, {\it Highly Degenerate Quasilinear Parabolic
Systems}. Ann. Scuola Sup. Pisa. Cl. Sci. {\bf 4}, 18, 135-166
(1991)

\item Bivar-Weinholtz and M. L. Lapidus, {\it Product Formula for
Resolvents of Normal Operators and the Modifed Feynman Integral}.
Proc. Amer. Math. Soc., Vol. 110, No. 2 (1990).

\item H. Br\'ezis, {\it Op\'erateurs maximaux monotones et
semi-groupes de contractions dans les espaces de Hilbert}, Math.
Studies 5, North Holland, Amsterdam, 1979

\item T. Diagana, {\it Sommes d'op\'erateurs et conjecture de
Kato-McIntosh}. C. R. Acad. Paris, t. 330, S\'erie I, p. 461-464
(2000).

\item G. Dore and A. Venni, {\it On the closedness of the Sum of two
closed operators}. Math. Z. 196, 189-201 (1987).

\item G. Fichera, {\it Linear Elliptic Differential Systems and
Eigenvalue Problems}. Lecture Notes in Math. 8, Springer-Verlag,
1965.

\item S. Kouachi, {\it Uniform Boundness and Global Existence of
Solutions for Reaction-Diffusion Systems with a Balance Law and a
Full Matrix of Diffusion}, Electon. J. Qual. Theory Differ. Equ.
(2001), No. 7, 9pp.

\item A. Pazy, {\it Semigroups of linear operators and applications to
partial differential equations}, Springer-Verlag, New York, 1983.

\item W. Rudin, {\it Functional analysis}, Tata McGraw-Hill, New Delhi,
1974.

\end{enumerate}

\vspace{1.2cm}

\noindent Toka Diagana, Howard University, Dept. of Mathematics,
2441 Sixth Street, N.W - Washington, DC 20059 - USA / E-mail:
tdiagana@howard.edu

\end{document}